\documentclass[12pt]{article}
\usepackage{theorem, amsfonts,amsmath,epsfig}

\nonstopmode
\def\CO{{\cal O}}

\def\CB{{\cal B}}
\def\CC{{\cal C}}

\def\CF{{\cal F}}
\def\CG{{\cal G}}
\def\CH{{\cal H}}

\def\CS{{\cal S}}

\def\pa{\partial}

\newcommand\COMP{\hbox{C\kern -.58em {\raise .54ex
\hbox{$\scriptscriptstyle |$}}\kern-.55em {\raise .53ex
\hbox{$\scriptscriptstyle |$}} }}

\newcommand\MM{\hbox{I\kern-.2em\hbox{M}}}
\newcommand\NN{\mathbb{N}}  
\newcommand\RR{\mathbb{R}}
\newcommand\sRR{{\it \hbox{I\kern-.2em\hbox{R}}}}
\newcommand\QQ{\hbox{I\kern-.53em\hbox{Q}}}
\newcommand\PP{\hbox{I\kern -.2em\hbox{P}}}
\newcommand\EE{\hbox{I\kern-.2em\hbox{E}}}
\newcommand\ZZ{{{\rm Z}\kern-.28em{\rm Z}}}
\newcommand\II{{{\rm I}\kern-.28em{\rm I}}}
\def\One{{ 1}\kern-.28em{\rm l}}

\def\proof{\noindent{\bf Proof:\ }}
\newcommand\qed{\hfill$\sqcap\kern-8.0pt\hbox{$\sqcup$}$\\}
\newcommand\be{\begin{equation}}
\newcommand\ee{\end{equation}}

\newtheorem{theorem}{Theorem}[section]
\newtheorem{proposition}[theorem]{Proposition}

\newtheorem{remark}[theorem]{Remark}

\newtheorem{corollary}[theorem]{Corollary}
\newtheorem{definition}[theorem]{Definition}

\title{Wiener chaos and the Cox--Ingersoll--Ross model}
\author{M. R. Grasselli$\ \!^*$ and T. R. Hurd\thanks{{Research supported by the
Natural Sciences and Engineering Research Council of Canada and Mathematics of
Information Technology and Complex Systems, Canada}} \\ Dept. of Mathematics and
Statistics\\ McMaster University\\Hamilton ON L8S 4K1\\Canada}
\begin{document}
\maketitle

\begin{abstract} In this we paper we recast the Cox--Ingersoll--Ross model
of interest rates into the chaotic representation recently introduced by
Hughston and Rafailidis. Beginning with the ``squared Gaussian
representation'' of the CIR model, we find a simple expression for the
fundamental random variable $X_\infty$. By use of techniques from the
theory of infinite dimensional Gaussian integration, we derive an
explicit formula for the $n$th term of the Wiener chaos expansion of the
CIR model, for $n=0,1,2,\dots$. We then derive a new expression for the
price of a zero coupon bond which reveals a connection between Gaussian
measures and Ricatti differential equations.  
\end{abstract}

\vspace{0.2in}

{\bf Key words:} Interest rate models, Wiener chaos, functional integrals, squared Gaussian models.

\newpage

\section{Introduction}

In the present paper we shall study the best known example of a term structure with
positive interest rates, namely the CIR model \cite{CoInRo85b}, in the context of the ``chaotic approach'' to interest rate 
dynamics introduced recently by Hughston and Rafailidis
\cite{HugRaf03} (see also \cite{BroHug03}). 
By an interest rate model we mean the specification of a spot rate process $r_t$ and of a market price of risk process
 $\lambda_t$ both under
the ``natural'' or physical measure in the economy $P$ . In the chaotic approach, the random nature of the model is assumed
to be given by a probability space $(\Omega,\CF,P)$ equipped with a Brownian
filtration $(\CF_t)_{0\leq t \leq \infty}$. The essence of the approach by Hughston and 
Rafailidis is the specification of the most general term structure with positive interest rates in terms of a 
single unconstrained random variable they denote $X_\infty$. They then apply a Wiener chaos expansion to  $X_\infty$, and 
interpret the resulting terms as building blocks for models of increasing complexity.

In one version of the CIR model, 
$r_t,\lambda_t$ are governed by the equations
\begin{equation}
\label{CIR2} dr_t=a(b-r_t)dt +c\sqrt{r_t} d\tilde W_t,\qquad
\|\lambda_t\|^2=\bar\lambda^2 r_t
\end{equation} for some positive constants
$a,b,c,\bar\lambda$ with $4ab> c^2$, where $\tilde W_t$ is a standard one dimensional $P$--Brownian
motion and $\lambda_t$ is a vector valued adapted process.  
By embedding this model inside a
general class of squared Gaussian models, we will be led to a natural
choice for the basic square integrable random variable
$X_\infty$ associated to it.  
As we shall then demonstrate, the resulting stochastic process
$X_t=E_t[X_\infty]$ admits an explicit chaos expansion, one which in general includes
terms of every chaotic order. 

A central idea in the present paper is the link between the Wiener chaos expansion and the theory of Gaussian functional 
integration, an
essential tool invented to study the mathematical structure of quantum
field theory. In fact, the mathematics underlying our example is a
consequence of certain basic results in that theory, and can be found in for example \cite{GliJaf81}. 

The organization of the paper is as follows. In section 2 we review the essential
ingredients for the construction of positive interest rates models in both
Flesaker--Hughston and the state price density approaches, and then compare these approaches to the
recently introduced chaotic representation of Hughston and Rafailidis \cite{HugRaf03}. 
We end the section by  describing the
structure of the Wiener--It\^{o} chaos expansion and show how it can be
expressed in terms of a certain generating functional acting on the space
$L^2(\RR^N)$.  

In section 3, we describe the squared Gaussian formulation of the  CIR
model and show the spot rate process can be explicitly
computed. Based on this representation, we state the form of the random
variable $X_\infty$, and give a proof that it lies in $L^2$. In section 4, we state the exponential 
quadratic formula which is the   main technical tool in this paper. It is a formula for the generating functional of  
random variables of the form $X=e^{-Y}$ for $Y$ lying in a 
general class of elements in the second chaos space $\CH_2$.   In section 5, we
compute the generating function for the random variables $X_t$ in the CIR model and, as the
main result of the paper, derive their chaos expansion.
In section 6,  we show that the usual CIR bond pricing formula has a natural
derivation within the chaotic framework.

Three appendices focus on the theory of Gaussian functional integration
and its relation to the Wiener chaos expansion.  Appendix A explores the white
noise calculus. Appendix B states and provides a proof of the generating
functional theorem. Appendix C provides a proof of the exponential quadratic formula.

\section{Positive Interest Rates} 
\label{positive}

\subsection{State price density and the potential approach}

Rather than focus on the spot rate process, one can model the system of bond prices directly. Let  $P_{tT}, 0\le t\le T$ denote
 the price at time
$t$ for a zero coupon bond which pays one unit of currency at its maturity $T$. Clearly $P_{tt}=1$ for
all $0\leq t < \infty$ and furthermore, positivity
of the interest rate is equivalent to having 
\begin{equation}
P_{ts}\leq P_{tu},
\label{bondpos}
\end{equation}
for all $0\leq t\leq u \leq s$.

A general way to model bond prices \cite{Rutk97,Roge97} is to write 
\begin{equation}
P_{tT}=\frac{E_t[V_T]}{V_t},
\label{FHkernel}
\end{equation}  
for a positive adapted continuous process $V_t$, called the {\em state price density}. Positivity of the interest rates is then equivalent to $V_t$ being a supermartingale. In order to match
the initial term structure, this supermartingale needs to be chosen so that $E[V_T]=P_{0T}$. If we further
impose that $P_{0T}\rightarrow 0$ as $T\rightarrow \infty$, then $V_t$ satisfies all the properties
of what is known in probability theory as a {\em potential} (namely, a positive  supermartingale
with expected value going to zero at infinity).

It follows from the Doob--Meyer decomposition that any continuous potential satisfying
\begin{equation}
E\left(\sup_{0\leq t \leq \infty} V^2_t\right) < \infty
\end{equation}
can be written as
\begin{equation}
V_t=E_t[A_\infty]-A_t,
\label{doob} 
\end{equation}
for a unique (up to indistinguishability) adapted continuous increasing process $A_t$ with $E(A^2_\infty)<\infty$.
Therefore, the model is completely specified by the process $A_t$, which can be freely chosen apart from 
the constraint that 
\begin{equation}
E\left[\frac{\partial A_T}{\partial T}\right]=-\frac{\partial P_{0T}}{\partial T}.
\end{equation}

\subsection{Related quantities and absence of arbitrage}

An earlier framework for positive interest rates was introduced by Flesaker and Hughston \cite[eq. (8)]{FleHug96}, who observed that any arbitrage free system of zero coupon bond prices has the form
\begin{equation}
P_{tT}=\frac{\int_T^\infty h_s M_{ts}ds}{\int_t^\infty h_s M_{ts}ds}, \quad \mbox{for } 
0\leq t \leq T<\infty.
\label{FHbonds}
\end{equation}
Here $h_T=-\frac{\partial P_{0T}}{\partial T}$ is a positive deterministic function obtained from the initial term structure and $M_{ts}$ is a 
family of strictly positive continuous martingales satisfying $M_{0s}=1$. Any such system of prices can be put into a potential form by
setting 
\begin{equation}
V_t=\int_t^\infty h_s M_{ts}ds.
\label{FHdensity}
\end{equation}
The converse result is less direct and was first established
by Jin and Glasserman \cite[lemma 1]{JinGla01}. 

These equivalent ways of
modelling positive interest rates can now be related to other standard financial objects. A particularly straightforward path is to follow 
\cite[proposition 1]{Rutk97}: given a strictly positive supermartingale
$V_t$, there exists a unique strictly positive (local) martingale $\Lambda_t$ such that the process
$B_t=\Lambda_t/V_t$ is strictly increasing. We identify $B_t$ with a riskless money market 
account initialized at $B_0=1$ and write it as
\begin{equation}
B_t=\exp(\int_0^t r_s ds),
\end{equation}
for an adapted process $r_s>0$, the short rate process. 

A sufficient condition for an arbitrage free bond price structure in the potential approach is to require that the local martingale
$\Lambda_t$ be in fact a martingale, since it can then be used as the density for an equivalent martingale measure. 
It is an interesting open question in the theory to isolate what conditions on the potential
$V_t$ would suffice for that.   

The formulation up to this point is quite general, in the sense that is does not make use of any particular
structure of the underlying filtration (other than the usual conditions). Let us now assume that ${\cal F}_t$ is
actually generated by an $N$-dimensional Brownian motion $W_t$. The market price of risk then arises as the adapted
vector valued process 
$\lambda_t$ such that 
\begin{equation}
d\Lambda_t=-\lambda_t\Lambda_t dW_t.
\end{equation}
It is also immediate to see that the state price density process is the solution to
\begin{equation}
dV_t=-r_tV_tdt - \lambda_tV_tdW_t,
\label{Vsde} 
\end{equation}
so that the specification of the process $V_t$ is enough to produce both the short rate
$r_t$ and the market price of risk $\lambda_t$.

It had already been remarked by Flesaker and Hughston \cite{FleHug96} that in the Brownian filtration
with finite time horizon any positive interest
rate model in their formulation corresponds to a model in the HJM family with positive 
instantaneous forward rates $f_{tT}$. The converse result that any interest rate model
in HJM form with positive instantaneous forward rates can be written in the Flesaker--Hughston form
was also obtained by Jin and Glasserman \cite[theorem 5]{JinGla01}. In order to prove this
result they found a rather technical necessary and sufficient condition for positivity in terms of
the volatility structure of the HJM form, confirming that the HJM formulation is not the most natural
one to investigate positive interest rates.

\subsection{The Chaotic Approach} 

We have seen in the potential approach that the fundamental 
ingredient to model the random behaviour of the interest rates is the increasing process $A_t$ in the decomposition
$V_t=E_t[A_\infty]-A_t$, whereas in the Flesaker--Hughston construction the corresponding role is played by the martingales $M_{ts}$.

In \cite{HugRaf03}, Hughston and Rafailidis introduced an elegant
construction of general positive interest rate models based on a Brownian filtration using simpler fundamentals.
Assuming that the state price density $V_t$ is a potential satisfying 
\begin{equation}
E\left[\int_0^\infty r_sV_s ds\right]<\infty,
\end{equation}
then integrating \eqref{Vsde} on the interval $(t,T)$, taking conditional expectations at time $t$
and the limit $T\rightarrow \infty$, one finds that 
\begin{equation}
V_t=E_t\left[\int_t^\infty r_sV_s ds\right].
\label{Vexp}
\end{equation}
Now let $\sigma_t$ be a vector valued process such that
\begin{equation}
\sigma_t\cdot \sigma_t=r_tV_t,
\label{sigma}
\end{equation}
and define the square integrable random variable 
\begin{equation}
\label{Xinf}
X_\infty=\int_0^\infty \sigma_s dW_s.
\end{equation}
It then follows from the It\^o isometry that
\begin{equation}
V_t=E_t[X_\infty^2]-E_t[X_\infty]^2,
\label{condvar}
\end{equation}
which is called the conditional variance representation of the state price density $V_t$. To obtain
the connection between this representation and the Flesaker--Hughston framework, observe that a direct
comparison between \eqref{Vexp} and \eqref{FHdensity} gives that 
\begin{equation}
h_s M_{ts}=E_t\left[\sigma_s\cdot \sigma_s\right].
\end{equation}
Similarly, by comparing the conditional variance representation \eqref{condvar} with the decomposition
\eqref{doob}, we see that 
\[E_t[X_\infty^2]-X_t^2=E_t[A_\infty]-A_t,\]
where $X_t=E_t[X_\infty]$. It follows from the uniqueness of the Doob-Meyer decomposition that 
\[A_t=[X,X]_t,\]
that is, the quadratic variation of the process $X_t$.

Conversely, given a zero-mean random variable $X_\infty\in L^2(\Omega,{\cal F},P)$, the representation
\eqref{condvar} defines a potential $V_t$, which can then be used as a state price density to obtain
a system of bond prices. The issue of absence of arbitrage can then be addressed in terms of necessary and sufficient
conditions on $X_\infty$, and is by and large an open question at this point. 
The construction in \cite{HugRaf03} flows in the opposite direction, in the
sense that the authors first enumerate a series of axioms to be satisfied by an arbitrage-free interest rate model
and then obtain a square integrable random variable $X_\infty$ corresponding to it.  

\subsection{Wiener chaos}

As Hughston and Rafailidis also observed, the $L^2$ condition on $X_\infty$ is necessary and 
sufficient for $X_\infty$ to have the type of orthogonal decomposition known as a Wiener chaos expansion \cite{Nual95}. 
They interpret the different orders of this decomposition as basic building blocks for models of increasing complexity. 

Let $W_t$ be an $N$--dimensional Brownian motion on the filtered
probability space
$(\Omega,\CF,\{\CF_t\}_{t\in \RR_+},P)$. We introduce a compact notation
\[\tau=(s,\mu)\in\Delta\doteq \RR_+\times\{1,\dots,N\}\]
and express
integrals as
\begin{eqnarray}
  \int_\Delta f(\tau) d\tau &\doteq&\sum_\mu \int^\infty_0 f(s,\mu)ds\nonumber\\
\int_\Delta f(\tau) dW_\tau &\doteq&\sum_\mu \int^\infty_0 f(s,\mu)dW^\mu_s
\end{eqnarray} 

For each $n\ge 0$, let 
\begin{equation}
\label{hermite}
H_n(x)= (-1)^n e^{x^2/2}\frac{d^n}{dx^n}e^{-x^2/2}
\end{equation}
 be the $n$th Hermite polynomial.  For $h\in   L^2(\Delta)$, let $\|h\|^2=\int_\Delta h(\tau)^2 d\tau$ and let
$W(h)$ be the Gaussian random variable
$\int_\Delta  h(\tau) dW_\tau$.
The spaces
\begin{eqnarray*}
\CH_n&\doteq& {\rm span}\{H_n(W(h))|h\in L^2(\Delta)\}, \quad n\ge 1, \\
\CH_0&\doteq& \mathbb{C}
\end{eqnarray*} 
 form an orthogonal decomposition of the space
$L^2(\Omega,\CF_\infty,P)$ of square integrable random variables:
\[L^2(\Omega,\CF_\infty,P)=\oplus_{n= 0}^\infty \CH_n\]

Each $\CH_n$ can be understood completely via the {\em isometries}
\[J_n:L^2(\Delta_n)\to\CH_n\]
given by
\begin{equation}
\label{In}f_n \mapsto J_n(f_n)=\int_{\Delta_n}
f_n(\tau_1,\dots,\tau_n)dW_{\tau_1}\dots dW_{\tau_n}
\end{equation} where $
\Delta_n\doteq
\{(\tau_1,\dots,\tau_n)|\tau_i=(s_i,\mu_i)\in\Delta,0\le s_1\le
s_2\le\dots\le s_n<\infty\}$.

With these ingredients, one is then led to the result that any $X\in
L^2(\Omega,\CF_\infty,P)$ can be represented as a {\em Wiener chaos
expansion}
\begin{equation} X=\sum_{n=0}^\infty J_n(f_n)
\end{equation} where the deterministic functions $f_n\in
L^2(\Delta_n)$ are uniquely determined by the random variable $X$.

A special example arises by noting  that for $h\in L^2(\Delta)$
\begin{equation}
\label{hermexp} n!J_n(h^{\otimes n})=\|h\|^n
H_n\left(\frac{W(h)}{\|h\|}\right)\end{equation} and furthermore
\begin{equation}
\exp\left[W(h)-\frac12\int h(\tau)^2d\tau\right]=\sum_{n=0}^\infty \frac{\|h\|^n}{n!}
H_n\left(\frac{W(h)}{\|h\|}\right)\end{equation} In the notation of quantum
field theory (see Appendix A), this example defines the Wick ordered
exponential and Wick powers
\begin{eqnarray} :\exp[W(h)]:\
&\doteq&\exp\left[W(h)-\frac12\int h(\tau)^2d\tau\right]\nonumber\\
:W(h)^n:&\doteq&n!J_n(h^{\otimes n}) 
\end{eqnarray}

Generating functionals provide one systematic approach to developing
explicit formulas for the terms of the chaos expansion in specific
examples.
\begin{theorem}
\label{genfnthm} For any random variable $X\in L^2(\Omega,\CF_\infty,P)$,
the generating functional $Z_X(h):L^2(\Delta)\to \mathbb{C}$ defined by
\begin{equation}
\label{generating} Z_X(h)\doteq  E\left[X \exp\left[W(h) -\frac12\int
h(\tau)^2 d\tau\right]\right]
\end{equation} is an entire analytic functional of $h\in L^2(\Delta)$ and
hence has an absolutely convergent expansion
\begin{equation}
 Z_X(h)=\sum_{n\ge 0}F_X^{(n)}(h)
\end{equation} where 
\begin{equation}
F_X^{(n)}(h)=\int_{\Delta_n}f_X^{(n)}(\tau_1,\dots,\tau_n)h(\tau_1)\dots
h(\tau_n)d\tau_1\dots d\tau_n.
\end{equation} The $n$-th Fr\'{e}chet derivative of $Z_X$ at
$h=0$, $ f_X^{(n)}(\tau_1,\dots,\tau_n)$, lies in $L^2(\Delta)$. 
Finally, the Wiener--It\^{o} chaos expansion of $X$ is
\begin{equation}
\label{Xexp} X=\sum_{n\ge 0} \int_{\Delta_n}
f_X^{(n)}(\tau_1,\dots,\tau_n)dW_{\tau_1} \dots dW_{\tau_n} 
\end{equation}
\end{theorem}

\proof See appendix B.

\qed

\medskip

\section{Squared Gaussian models}

A number of authors  \cite{Magh96},\cite{Jams95}, \cite{Roge97} have observed that the CIR
model \cite{CoInRo85b} with an integer constraint $N\doteq
\frac{4ab}{c^2}\in\NN_+\setminus\{0,1\}$ lies in the class of so--called
squared Gaussian models. By introducing an
$\RR^N$--valued Ornstein--Uhlenbeck process $R_t$, governed by the
stochastic differential equation
\begin{equation}
 dR_t=-\frac{a}2 R_t\ dt + \frac{c}2 dW_t
 \label{sqgauss}
\end{equation} where $W_t$ is $N$--dimensional Brownian motion,  the It\^{o} formula together with L\'evy's criterion for Brownian motion shows that the square $r_t=R^\dagger_tR_t$
satisfies (\ref{CIR2})  where $\tilde W_t=\int^t_0 (R^\dagger_t R_t)^{-1/2}R_t\cdot
dW_t$  is itself a one--dimensional Brownian motion. Note that here, and for the remainder of the paper, we suppress vector indices by adopting a matrix multiplication convention, including $\dagger$ for transpose, in which for example $ R^\dagger_t\gamma(t)dW_t$ denotes $\sum_{\mu,\nu=1}^N R^\mu_t\gamma^{\mu\nu}(t)dW^\nu_t$. 

We focus on a general
family of interest models  which includes this example, the so-called extended CIR model and more. Note that we work always in the physical measure and thus to specify  the term structure model one needs to
determine  the market price of risk vector $\lambda_t$ as well as the spot rate process $r_t$.
\begin{definition} A pair $(r_t,\lambda_t)$ of  $(\Omega,\CF,\CF_t,P)$  processes is
called an
$N$--dimensional squared Gaussian model of interest rates ($N\ge 2$) if  there is an $\RR^N$--valued
Ornstein--Uhlenbeck process  such that $r_t=R^\dagger_t R_t$ and $\lambda_t=\bar\lambda R_t$.  $R_t$ satisfies
\begin{equation}
\label{OUeqn}
  dR_t=\alpha(t)(\bar R(t)-R_t)dt+\gamma(t) dW_t,\quad R|_{t=0}=R_0
\end{equation}
where $\alpha, \gamma, \bar\lambda$ are symmetric matrix valued and $\bar
R$ vector valued deterministic measurable functions on $\RR_+$.  $W$ is standard $N$--dimensional Brownian motion. In
addition we impose boundedness conditions:
\begin{itemize}
  \item  there is some constant $M>0$ such that
$ \alpha(t)\ge M$ and $|\bar\lambda(t)|^{-2} \ge M$ for all $t$
 \end{itemize}
\end{definition}

The exact solution of (\ref{OUeqn}) is easily seen to be 
\begin{equation}
\label{Rform} R_t=\tilde R(t)+\int  K(t,t_1)(\gamma dW)_{t_1}
\end{equation} where
\begin{equation}
\tilde R(t)=K(t,0) R_0+\int K(t,t_1)\alpha(t_1)\bar
R(t_1) dt_1 
\label{tildeRform}
\end{equation}
and $ K(t,s), t\ge s$ is the matrix valued solution of 
\begin{equation}
\label{Kgen }
\left\{\begin{array}{ll }
  dK(t,s)/dt=-\alpha(t)K(t,s)    &0\le s\le t    \\
  K(t,t)=I    &   0\le t
\end{array}\right.
\end{equation} which generates the Ornstein--Uhlenbeck semigroup.

By (\ref{Vsde}), the state price density process is
\begin{equation} V_t=\exp\Bigl[-\int^t_0\biggl(R_s^\dagger\left(1+\frac{\bar\lambda^2}{2}\right) R_s
ds+ R_s^\dagger\bar\lambda dW_s\biggr)\Bigr]
\end{equation}  We thus have a natural candidate for the random variable
$X_\infty$:
\begin{equation}\label{Xinfty} X_\infty=\int^\infty_0 \sigma_t^\dagger\ dW_t
\end{equation}
where the $\RR^N$--valued process 
\begin{equation} \sigma_t\doteq
\exp\biggl[-\int^t_0
\biggl(R_s^\dagger\left(\frac{1}{2}+\frac{\bar\lambda^2}{4}\right) R_s
ds+\frac12 R_s^\dagger\bar\lambda dW_s\biggr)\biggr] R_t
\end{equation} 
is the natural solution of $\sigma^\dagger _t\sigma_t=r_tV_t $.

Before proceeding to analyse $X_\infty$ in detail  we show that
$X_\infty$ is square integrable.

\begin{proposition} $ \lim_{T\to\infty}E[X_T^2] =1$
\end{proposition}

\proof By the It\^{o} isometry and Fubini's theorem for the It\^{o}
integral,
\begin{eqnarray} E[X_T^2] &=&E\left[\int^T_0
\sigma^\dagger _t\sigma_t dt\right] 
\end{eqnarray} From
$\sigma^\dagger _t\sigma_t dt =r_tV_t dt=-dV_t- V_t \lambda^\dagger dW_t
$ it follows that
$ E[X_T^2] = 
E\left[1-V_T \right]
$. We observe that this is true for any pair $X_t$, $V_t$ defined using the formal relations
\eqref{Vsde},\eqref{sigma},\eqref{Xinf} of the chaotic approach. All we need to prove now is that the 
supermartingale $V_t$ is indeed a potential when $(r_t,\lambda_t)$ is a squared Gaussian model.

For any $0<\epsilon<2M$, $V_T=e^{-Y_1-Y_2}$ where
$Y_1=\int^T_0R_t^\dagger \bigl(1-\frac{\epsilon\bar\lambda^2}2\bigr)R_tdt$
is positive and
$$Y_2=\frac{1}{1+\epsilon}\int^T_0\left[\frac12 R_t^\dagger (1+\epsilon)^2\bar\lambda^2R_tdt
+(1+\epsilon) R_t^\dagger \bar\lambda^\dagger dW_t\right]$$ 
By the H\"older inequality
$$E[V_T]\le\left(E[e^{(1+1/\epsilon)Y_1}]\right)^{\epsilon/(1+\epsilon)}\left(E[e^{-(1+\epsilon)Y_2}]\right)^{1/(1+\epsilon)}$$
with the second factor equal to $1$, since $e^{-(1+\epsilon)Y_2}$ is an exponential martingale. Now $Y_1$ is a positive random variable for which a direct
computation shows
\begin{eqnarray}
\mbox{mean}(Y_1)&=&C_1NT\left(1+\CO(T)\right)\\
\mbox{var}(Y_1)&=&C_2NT\left(1+\CO(T)\right)
\end{eqnarray} for positive constants $C_1,C_2$. An easy application of
Chebyshev's inequality
\begin{equation} \mbox{Prob}\left(Y_1\le \frac{C_1NT}2\right)\le
\CO(\frac1{NT})
\end{equation} then implies that $\lim_{T\to\infty}E[V_T]=0$.

\qed

\section{Exponentiated second chaos}
The chaos expansion we seek for the CIR model will be derived from a closed
formula for expectations  of $e^{-Y}$ for elements 
\begin{equation}
\label{Ydef}
Y=A+  \int_\Delta B(\tau_1)dW_{\tau_1}
+\int_{\Delta_2}C(\tau_1,\tau_2)dW_{\tau_1}dW_{\tau_2}
\end{equation}
in a certain subset $\CC^{+}\subset \CH_{\le 2}\doteq
\CH_0\oplus\CH_1\oplus\CH_2$.  In the integrals above, recall that compact notation using $\tau$'s carries a summation over vector indices as well as 
integration over time. The formula we present is well known in the theory of Gaussian functional integration \cite[Chapter 9]{GliJaf81}. In probability theory, this result gives the Laplace transform of a general class of quadratic functionals of Brownian motion. Many special cases of this result have been studied in probability theory, see for example \cite[Chapter 2]{Yor92a} and the references contained therein.

If in (\ref{Ydef})  we define $C(\tau_1,\tau_2)=C(\tau_2,\tau_1)$ when $\tau_1>\tau_2$, then $C$ is the kernel of a symmetric integral operator on $L^2(\Delta)$:
\begin{equation}
[Cf](\tau)=\int^\infty_0 C(\tau,\tau_1)f(\tau_1)d\tau_1 
\end{equation}
Recall that Hilbert-Schmidt operators on $L^2(\Delta)$ are finite norm operators under
the norm:
\[ \|C\|^2_{HS}=\int_{\Delta^2} C(\tau_1,\tau_2)^2d\tau_1d\tau_2
\]
We say that $Y\in\CH_{\le 2}$ is in $\CC^{+}$ if $C$ is the
kernel of a symmetric Hilbert--Schmidt operator on $L^2(\Delta)$
such that $(1+C)$ has positive spectrum.
\begin{proposition}
\label{measchange}
  Let $Y\in\CC^+ $.  Then
\begin{eqnarray}
\label{change} E [e^{-Y}]&=&\left[{\rm det}_2 (1+
C)\right]^{-1/2}\nonumber\\ &&\hspace{-.5in}
\exp\left[-A+ \frac12\int_{\Delta_2}
B(\tau_1)(1+C)^{-1}(\tau_1,\tau_2)B(\tau_2)
d\tau_1d\tau_2\right]
\end{eqnarray}
\end{proposition}

\begin{remark} The Carleman--Fredholm
determinant is defined as the extension of the formula
\begin{equation}
{\rm det}_2 (1+ C) = {\rm det} (1+ C)\exp[{-{\rm Tr}(C)}]
\end{equation} from finite rank operators to bounded 
Hilbert--Schmidt operators; the operator kernel
$(1+C)^{-1}(\tau_1,\tau_2)$ is also the natural extension from the finite rank case.   \end{remark}

\proof See Appendix C.

\qed

Using this proposition, it is possible to deduce the chaos expansion of the random variables 
$e^{-Y}, Y\in \CC^{+}$, a result known in quantum field theory as Wick's Theorem:

\begin{corollary}
\label{Feynmancor}
If  $Y=\int_{\Delta_2}C(\tau_1,\tau_2)dW_{\tau_1}dW_{\tau_2}\in \CC^{+} $, then the random variable $X=e^{-Y}$
has Wiener chaos coefficient functions
\begin{equation}
f_n(\tau_1,\dots,\tau_n)=\left\{
\begin{array}{ll} K\sum_{G\in\CG_n}\prod_{g\in G}
[ C(1+C)^{-1}](\tau_{g_1},\tau_{g_2})&n\
\mbox{even}\nonumber\\ 0&n\ \mbox{odd}\nonumber\\
\end{array}\right.
\end{equation} where $K=\left[\rm{det}_2 (1+C)\right]^{-1/2}$
and for $n$ even,  $\CG_n$ is the set of {\em Feynman graphs} on the
$n$ marked points
$\{\tau_1,\dots,\tau_n\}$. Each Feynman graph $G$ is a disjoint union of unordered
pairs
$g=(\tau_{g_1},\tau_{g_2})$ with $\cup_{g\in G}\
g=\{\tau_1,\dots,\tau_n\}$.

\end{corollary}

\proof The generating functional for $X=e^{-Y}$ is
\begin{eqnarray*}
Z_X(h)&=&E\left[X\exp\left(\int h(\tau)dW_\tau - \frac{1}{2}\int h(\tau)^2d\tau\right)\right]\\
&=& E\left[\exp\left(\int h(\tau)dW_\tau- \frac{1}{2}\int h(\tau)^2d\tau-\int_{\Delta_2}
C(\tau_1,\tau_2)dW_{\tau_1}dW_{\tau_2}\right)\right],
\end{eqnarray*}
so we can use Proposition \ref{measchange} with $A=\frac{1}{2}\int h(\tau)^2 d\tau$ and $ B(\tau)=-h(\tau)$, which yields
\begin{eqnarray}
\label{Z_X }
Z_X(h)&=&\mbox{det}_2 (1+C)^{-1/2}\nonumber\\
&&\hspace{-.3in}
\exp\left[-\frac12\int_{\Delta^2} h^\dagger(\tau_1)[\delta(\tau_1,\tau_2)-(1+C)^{-1}(\tau_1,\tau_2)]h(\tau_2)d\tau_1d\tau_2
\right].
\end{eqnarray}
Using the last part of Theorem \ref{genfnthm}, the result comes by evaluating the 
$n$th Fr\'echet derivative at $h=0$, or equivalently by expanding the 
exponential and symmetrizing over the points $\tau_1,\dots,\tau_{n}$ in the $n/2$th term.

\qed

\section{The chaotic expansion for squared Gaussian models}

We now derive the chaos expansion for the squared Gaussian model defined
by (\ref{OUeqn}).  In view of (\ref{Xinfty})  it will be enough to find the chaos expansion for
$\sigma_T^\mu$, $T<\infty$. 
We start by finding its generating functional $Z_{\sigma^\mu_T}$. 
For $h,k \in L^2(\Delta)$, define the auxiliary
functional $ Z(h,k) = E\left[e^{-Y_T}\right]$ with
\begin{eqnarray}Y_T&=&
\int^T_0R^\dagger_t\left(\frac{1}{2}+\frac{\bar\lambda^2}{4}\right)R_tdt +
\frac{1}{2}\int_0^TR^\dagger_t\bar\lambda dW_t - \int^T_0 h^\dagger(t) dW_t\nonumber\\ 
&&\hspace{.6in}+\frac12 \int^T_0 h^\dagger(t)h(t) dt -\int^T_0k^\dagger(t)R_t dt.
\label{Zhk}
\end{eqnarray}

\begin{proposition}
$Z(h,k)$ is an entire analytic functional on $L^2(\Delta)\times L^2(\Delta)$. 
Moreover
\begin{equation}\label{Zsigma}
\lim_{t\to T^-}\frac{\delta Z(h,k)}{\delta
k^\mu(t)}\Big|_{k=0}=Z_{\sigma^\mu_T}(h)
\end{equation} where $Z_{\sigma^\mu_T}(h)$ is defined by (\ref{generating})
with
$X=\sigma^\mu_T, \mu=1,\dots, N$.
\end{proposition}

\proof Analyticity in $(h,k)$ follows by repeating the argument given in
Appendix B. By the definition of Fr\'echet differentiation and continuity
of the
$t\to T^-$ limit, (\ref{Zsigma}) follows.

\qed

We want to use Proposition \ref{measchange} in order to compute $Z(h,k)$. Substitution of 
(\ref{Rform}), into the first term of (\ref{Zhk}) leads to 
\begin{eqnarray*} 
\int^T_0 R^\dagger_t\left(\frac12+\frac{\bar\lambda^2}{4}\right)R_tdt &=&
\int^T_0 \tilde R^\dagger(t)\left(\frac12+\frac{\bar\lambda^2}{4}\right)\tilde R(t)dt \\
&& +\int_0^T\left[\int_0^T\tilde R^\dagger(s)\left(1+\frac{\bar\lambda^2}{2}\right) K_T(s,t)ds\right]\gamma(t)dW_{t}\\
&& \hspace{-1.1in}+\int_{\Delta_2}\gamma(t_1)\left[\int_0^T K^\dagger_T(t_1,s)\left(1+\frac{\bar\lambda^2}{2}\right) K_T(s,t_2)ds\right]\gamma(t_2)dW_{t_1}dW_{t_2}\\
&&\hspace{-0.9in}+\int^T_0{\rm tr}\left\{\gamma(t)\left[ \int_0^T K_T^\dagger(t,s) \left(\frac12+\frac{\bar\lambda^2}{4}\right)K_T(s,t)ds\right]\gamma(t)\right\}dt,  \\
\end{eqnarray*} 
where we define $K_T(t_1,t_2)=\One(t_1\le T) K(t_1,t_2)$. For the second and the last terms of \eqref{Zhk} we have
\begin{eqnarray*}
\frac12\int_0^TR^\dagger_t\bar\lambda dW_t &=& \frac12\int_0^T\tilde{R}^\dagger(t)\bar\lambda dW_t 
+\frac12\int_0^T \mbox{ tr}\left[\int_0^T \gamma(s) K_T^\dagger(s,t)\bar\lambda ds\right]dt  \\ 
&&+\frac12\int_{\Delta_2}\left(\gamma(t_1) K_T^\dagger(t_1,t_2)\bar\lambda+\bar\lambda K_T(t_1,t_2)\gamma(t_2)\right) dW_{t_1}dW_{t_2},\\
\int^T_0k^\dagger(t)R_t dt &=& \int^T_0 k^\dagger(t) \tilde R(t) dt+\int_0^T \left( \int_0^T k^\dagger (s) K_T(s,t) ds\right) \gamma(t) dW_t.
\end{eqnarray*}
Thus  the exponent
$Y_T$ appearing in (\ref{Zhk}) has the form of (\ref{Ydef}) with
\begin{eqnarray*}
A_T&=&\int^T_0\left[\tilde R^\dagger(t) \left(\frac{1}{2}+\frac{\bar\lambda^2}{4}\right)\tilde R(t) +\frac12 h^\dagger(t)h(t) - k^\dagger(t)\tilde R(t)\right]dt  \\
&&+\int_0^T\mbox{tr}\left\{\gamma(t)\left[\int_0^T K_T^\dagger(t,s)\left(\frac{1}{2}+\frac{\bar\lambda^2}{4}\right)K_T(s,t)ds\right]\gamma(t)
\right\}dt \\
&&\frac12\int_0^T \mbox{tr}\left[\int_0^T \gamma(s) K^\dagger_T(s,t)  \bar\lambda ds \right] dt, \\
B_T(t)&=&-h(t)-\gamma(t)\int_0^T K^\dagger_T(t,s)k(s)ds + \frac12 \bar\lambda \tilde R(t)\\
&&+ \gamma(t)\int_0^T  K^\dagger_T(t,s)\left(1+\frac{\bar\lambda^2}{2}\right) 
\tilde R(s)ds  \\
C_T(t_1,t_2)&=& \gamma(t_1)\left[\int_0^T K_T^\dagger(t_1,s)\left(1+\frac{\bar\lambda^2}{2}\right) K_T(s,t_2)ds\right]\gamma(t_2) \\
&& +\frac{1}{2}\left[\gamma(t_1) K^\dagger_T(t_1,t_2)\bar\lambda+\bar\lambda K_T(t_1,t_2)\gamma(t_2)\right]
\end{eqnarray*}

It is clear that the operator $C_T$ has Hilbert-Schmidt norm $\|C_T \|_{HS}^2=\CO(T)$. Moreover, 
if we denote by $\gamma K^\dagger_T \left(1+\frac{\bar\lambda^2}{2}\right) K_T \gamma$, 
$\gamma K^\dagger_T\bar\lambda$ and $\bar\lambda K_T\gamma$ the operators whose kernels appear in the expression
above, then $C_T$ can be written as
\begin{eqnarray}
C_T &=& \gamma K_T^\dagger K_T \gamma + \frac12(\gamma K^\dagger_T\bar\lambda+1)(\bar\lambda K_T\gamma+1)-\frac{1}{2}
\end{eqnarray}
from which we see that $(1+C_T)$ is positive. Therefore, we can use Proposition \ref{measchange} for 
$E[e^{-Y_T}]$, leading to a general formula for the generating functional
$Z(h,k)$:

\begin{eqnarray}
Z(h,k)&=& \mbox{det}_2 (1+C_T)^{-1/2}\exp{\left(-\frac{1}{2}\mbox{tr} C_T\right)} \nonumber \\
&&\hspace{-0.5in}\times\exp\left\{-\int_0^T \left[ \tilde 
R^\dagger(t)\left(\frac{1}{2}+\frac{\bar\lambda^2}{4}\right) \tilde R(t) + \frac{1}{2} h^\dagger(t)h(t)- k^\dagger(t) 
\tilde R(t)  \right]dt \right\}  \nonumber\\
&& \hspace{-.5in}\times \exp\left\{\frac{1}{2}\int_{\Delta_2} \left[h^\dagger(t_1) + \int_0^T k^\dagger(s) K_T(s,t_1) \gamma(t_1)ds
-\frac12\tilde R^\dagger(t_1)\bar\lambda \right.\right.  \nonumber\\
&&\left. -\int_0^T  \tilde R^\dagger(s) \left(1+\frac{\bar\lambda^2}{2}\right)K_T(s,t_1)\gamma(t_1)ds\right](1+C_T)^{-1}(t_1,t_2) \nonumber\\
&&\times  \left[h(t_2)  + \gamma(t_2) \int_0^T K^\dagger_T(t_2,s) k_s ds -\frac12 \bar\lambda\tilde R(t_2) \right.  \nonumber\\
&&\left.\left.- \gamma(t_2) \int_0^T  K^\dagger_T(t_2,s)\left(1+\frac{\bar\lambda^2}{2}\right)\tilde R(s) ds\right]dt_1dt_2 \right\}
\label{mess1}
\end{eqnarray}

Differentiation once with respect to $k$  then yields 
\begin{eqnarray}
Z_{\sigma_T}(h)&=&M_T\exp\left\{-\int^T_0\left[ \tilde R^\dagger(t)\left(\frac{1}{2}+\frac{\bar\lambda^2}{4}\right) \tilde R(t)+
\frac12 h^\dagger(t) h(t) \right]dt\right\}  \nonumber\\
&& \hspace{-0.3in}\times\left\{-\tilde R+ K^{}_T\gamma(1+C_T)^{-1}\left[h-\frac{\bar\lambda\tilde R}{2}-\gamma K^\dagger_T\left(1+\frac{\bar\lambda^2}{2}\right)
\tilde R\right]\right\}(T)  \nonumber\\
&&\hspace{-0.3in}\times\exp\left\{\frac{1}{2}\int_{\Delta_2}\left[h^\dagger -\frac{\tilde R^\dagger\bar\lambda}{2}
-\tilde R^\dagger\left(1+\frac{\bar\lambda^2}{2}\right) K^{}_T\gamma\right](t_1)
(1+C_T)^{-1}(t_1,t_2)\right.  \nonumber\\
&&\hspace{0.1in}\left.\times\left[h-\frac{\bar\lambda\tilde R}{2}-\gamma K_T^\dagger\left(1+\frac{\bar\lambda^2}{2}\right)
 \tilde R\right](t_2)dt_1dt_2\right\} 
\label{mess2}
\end{eqnarray} 
where
\begin{equation} M_T=e^{-\frac12{\rm tr}C_T}({\rm det}_2(1+C_T))^{-1/2}=({\rm det}(1+C_T))^{-1/2}.
\end{equation}
By comparing \eqref{mess1} and \eqref{mess2}, the reader can observe our use of an operator notation which suppresses
some time integrals, for example in the very last term
\[\left[\gamma K_T^\dagger\left(1+\frac{\bar\lambda^2}{2}\right) \tilde R\right](t)\doteq
\gamma(t) \int_0^T  K^\dagger_T(t,s)\left(1+\frac{\bar\lambda^2}{2}\right)\tilde R(s) ds,\]
and similarly for other terms.

These formulas simplify considerably if the function
$\tilde R$ vanishes, which is true in the simple CIR model of \eqref{sqgauss} when $r_0=0$. In this case
we have $\alpha(t)=a/2$ and $\gamma(t)=c/2$ (here certain scalars are to be understood as multiples of the identity matrix), so that $K_T(s,t)=e^{-a(s-t)/2}\One (t\leq s\leq T)$ and 
\begin{eqnarray}
\!\!\! C_T(t_1,t_2)
&=&\frac{c^2}{4a}\left(1+\frac{\bar\lambda^2}{2}\right)\left[e^{-\frac{a}{2}|t_1-t_2|}-e^{\frac{a}2(t_1+t_2-2T)}\right]+\frac{c}{2}\bar\lambda 
e^{-\frac{a}{2}|t_1-t_2|}
\label{CIRc}
\end{eqnarray}
Moreover, the previous expression for $Z_{\sigma_T}(h)$ reduces to
\begin{eqnarray*}
Z_{\sigma_T}(h)&=&M_T\left[K_T\gamma(1+C_T)^{-1}h \right](T)  \\
&& \hspace{-0.4in} \exp\left[-\frac12\int^T_0 h^\dagger(t) h(t) dt + \frac{1}{2}\int_{\Delta_2}h^\dagger(t_1)
(1+C_T)^{-1}(t_1,t_2)h(t_2)dt_1dt_2\right] .
\end{eqnarray*}

We can then easily evaluate the $n$th Fr\'echet derivative of $Z_{\sigma_T}$ at $h=0$ as in the
 proof of Corollary \ref{Feynmancor} and determine the following partly explicit form for the
$n$th term of the chaos expansion. 

\begin{theorem} 
 The 
$n$th term of the chaos expansion of $\sigma_T$ for the
CIR model with initial condition $r_0=0$ is zero for $n$ even. For $n$ odd, the kernel of 
the expansion is the function
$f^{(n)}_T(\cdot):\Delta_n\to\RR$
\begin{equation} f_T(t_1,\dots,t_n)=M_T
\sum_{G\in\CG_n^*}\prod_{g\in G}L(g)
\end{equation} 
where
\begin{equation} L(g)=
\left\{\begin{array} {ll}
[C_T(1+C_T)^{-1}](t_{g_1},t_{g_2})
& T\notin g        
\\\\
(K_T\gamma(1+C_T)^{-1})(T,t_{g_2}) & T\in g
\end{array}  \right.
\end{equation}
Here,
$\CG_n^*$ is the set of {\em Feynman graphs}, each Feynman graph $G$ being a partition of $\{t_1,\dots,t_n,T\}$ into pairs $g=(t_{g_1},t_{g_2})$.
\end{theorem}

The chaos expansion for
$X_\infty$ itself is exactly the same, except that the variable $T$ is treated as an additional It\^o integration variable. The explicit expansion up to fourth order is:
\begin{eqnarray}
X_\infty&=&\int_{\Delta_2}M_T[K_T\gamma(1+C_T)^{-1}](T,t_1)dW_{t_1}dW_{T}
\nonumber\\
&&\hspace{-.5in}+\int_{\Delta_4} M_T[K_T\gamma(1+C_T)^{-1}](T,t_3)[C_T(1+C_T)^{-1}](t_1,t_2)
dW_{t_1}dW_{t_2}dW_{t_3}dW_{T}
\nonumber\\
&&\hspace{-.5in}+\int_{\Delta_4} M_T[K_T\gamma(1+C_T)^{-1}](T,t_2)
[C_T(1+C_T)^{-1}](t_1,t_3)dW_{t_1}dW_{t_2}dW_{t_3}dW_{T}
\nonumber\\
&&\hspace{-.5in}+\int_{\Delta_4} M_T[K_T\gamma(1+C_T)^{-1}](T,t_1)
[C_T(1+C_T)^{-1}](t_2,t_3)dW_{t_1}dW_{t_2}dW_{t_3}dW_{T} \nonumber\\ &&\cdots
\end{eqnarray}

\section{Bond pricing formula}

In this section we give a derivation of the price of a zero coupon bond
in the CIR model. Recall from section \ref{positive} that these are given by 
\begin{equation} P_{tT}=E_t[V_t^{-1} V_T].
\end{equation}
To keep things ``as
simple as possible, but not any simpler'', we take
 $\bar\lambda=0$ so $V_t=\exp[-\int^T_0 r_s ds]$, or in terms of the squared Gaussian formulation,
\begin{equation}
V_t=\exp\left[-\int_0^t R^\dagger_s R_s ds\right].
\end{equation} 
As we have seen in the previous section, for $t\le s \le T$ 
\begin{equation}
R_s^\mu=K_T(s,t)R_t^\mu+\frac{c}{2}\int_t^sK_T(s,s_1)dW_{s_1}^\mu,
\end{equation}
hence $-\log[V_t^{-1}V_T]=\sum_\mu \int^T_t(R^\mu_s)^2ds$ can be
written as
\begin{eqnarray}
\sum_\mu\left[\frac{4}{c^2}(R^\mu_t)^2 C_T(t,t) +  \frac{4R^\mu_t}{c}
\int^T_t C_T(t,s) dW^\mu_s +2\int^T_t\int^{s_2}_t C_T(s_1,s_2)dW^\mu_{s_1}dW^\mu_{s_2}\right]&&\nonumber\\
&&\hspace{-3in}+N \int^T_t C_T(s,s)ds 
\end{eqnarray} 
where $C_T(s_1,s_2)$ is given by \eqref{CIRc} with $\bar\lambda=0$. 

Taking the conditional expectation of $V_t^{-1} V_T$ by use of
Proposition \ref{measchange} leads to the desired formula
\begin{eqnarray}
P_{tT}&=&\left[\det(1+2C_T)\right]^{-N/2}\\&&
\prod_\mu\exp\left[-\frac{4}{c^2} (R^\mu_t)^2
\left(C_T (1+2C_T)^{-1}\right)(t,t)\right]
\end{eqnarray}
Thus $P_{tT}$ has the exponential affine form
$\exp[-\beta(t,T)r_t-\alpha(t,T)]$ with
\begin{eqnarray}
\beta(t,T)&=&\frac{4}{c^2}[C_T(1+2C_T)^{-1}](t,t)\nonumber\\
\alpha(t,T)&=&\frac{N}2\log\bigl[\det(1+2C_T)\bigr]\label{betaalpha}
\end{eqnarray}

The known formula has the same form, with 
\begin{eqnarray}
\beta(t,T)&=&\frac{2(e^{\rho(T-t)}-1)}{(\rho+a)(e^{\rho(T-t)}-1)+2\rho},
\quad\rho^2=a^2+2c^2 \nonumber\\
\alpha(t,T)&=& -\frac{2ab}{\rho^2} \log\left[\frac{2\rho e^{(a+\rho)(T-t)/2}}{(\rho+a)(e^{\rho(T-t)}-1)+2\rho}\right]
\end{eqnarray}
which can be derived as solutions of the pair of Ricatti ordinary
differential equations
\begin{eqnarray}
\frac{\pa \beta}{\pa  t}&=&\frac{c^2\beta^2}2+a\beta-1\nonumber\\
\frac{\pa \alpha}{\pa t}&=&  -ab\beta 
\end{eqnarray}

One can demonstrate using power series expansions that (\ref{betaalpha}) do in fact solve the
Ricatti equations and hence agree with the usual formula. This example
points to the rather subtle general relationship between kernels such as
$(1+2C_T)^{-1}$ and solutions of Ricatti equations deserving of further
study.

\section{Discussion}
We have shown how the CIR model, at least in integer dimensions, can be viewed within the chaos framework of Hughston and Rafailidis as arising from a somewhat special random variable $X_\infty$. This random variable can be understood as derived from exponentiated second chaos random variables $e^{-Y},Y\in \CC^+)$. Such exponentiated $\CC^+$ variables form a rich and natural family which is likely to include many more candidates for applicable interest rate models.  Although their analytic properties are complicated, there do exist approximation schemes which can in principle be the basis for numerical methods. 

On the theoretical side, this family is distinguished by its natural invariance properties. Most notably, as will be investigated elsewhere, it is invariant under conditional $\CF_t$--expectations: $\log E_t[e^{-Y}]\in \CC^+$ whenever $Y\in \CC^+$. Note in particular that this implies that these can be used as  the Radon--Nikodym derivatives of measure changes which generalize the Girsanov transform, and which can greatly enrich the tools  applicable in finance.

This application of the chaos expansion to squared Gaussian models also illustrates
a deep connection between methods developed for quantum field theory and the
methods of Malliavin calculus. Many of the very rich
analytic properties of this example reflect well known techniques widely used in mathematical physics.

To conclude, even if the representation of the CIR
model and its generalizations we present is not simple, it does show the way to the use of a powerful set of analytical techniques which might prove to be very useful in financial modelling.  

\appendix
\section{White noise calculus}

Here we describe the white noise calculus, which can be regarded as a reformulation of
the calculus of Wiener measure into concepts familiar to practitioners in quantum field
theory as Gaussian functional integration. We follow the discussion of {\O}ksendal
\cite{Okse96}.

Let $\CS$ be the Schwartz space of smooth functions on
$\RR_+$ of rapid decrease, and $\CS'$ its topological dual, the space of tempered
distributions on $\RR_+$. If $\phi\in\CS'$ and $f\in\CS$, we write
\[\phi(f)\doteq \langle\phi,f\rangle\]
for their canonical pairing. We also use the formal notation
\[\phi(f)=\int_{\RR_+} f(s)\phi_s ds\]
to represent the `smearing' of the distribution $\phi$ over the
test function $f$. It acquires a rigorous meaning, however, in
the cases where $\phi$ is itself a function on $\RR_+$ for which
the pointwise product $\phi(s)f(s)$ is integrable for all
$f\in\CS$.

Define now the functional
\begin{equation}
S\{f\}\doteq e^{-\frac{1}{2}\|f\|^2}=e^{-\frac{1}{2}\langle f,f
\rangle_{L^2}},
\end{equation}
where $\langle \cdot,\cdot \rangle_{L^2}$ denotes the real-valued
inner product in $L^2(\RR_+)$. Observe that this functional
satisfies
\begin{enumerate}
\item (continuity) $f_n\rightarrow f$ in $\CS$ implies that
$S\{f_n\}\rightarrow S\{f\}$ ;

\item (positive definiteness) $\sum_{i,j=1}^N
\overline{c_i}c_jS\{f_i-f_j\}\geq 0$ for all $f_i\in\CS$ and all
$c_i\in\mathbb{C}$  ;

\item (normalization) $S\{0\}=1$.
\end{enumerate}

It follows from the Bochner-Minlos theorem that there exists a
unique Borel probability measure $\mu$ on $\CS^\prime$  such that
$S\{f\}$ corresponds to a {\em moment generating functional}, that
is, for all $f\in\CS$
\begin{equation}
\int_{\CS'} e^{i\phi(f)} d\mu(\phi) =S\{f\}=e^{-\|f\|^2/2}.
\end{equation}
The measure space $(\CS',\CB,\mu)$, where $\CB$ is the Borel
$\sigma$-algebra of $\CS'$, is called the {\em white noise
probability measure}.

In Euclidean quantum field theory, a given Borel measure $\mu$ on
$\CS^\prime$ characterizes the family of `fields'
$\phi\in\CS^\prime$ through the properties of the random variables
$\phi(f):\CS^\prime\rightarrow\RR$ obtained for each $f\in\CS$.
One tries to construct measures $\mu$ so that the generating
functional $S\{f\}$ satisfy the so-called Osterwalder-Schrader
axioms, in order to guarantee that the fields $\phi$ have certain
required physical properties.   The family of Euclidean free fields is obtained when
\[S_C\{f\}\doteq e^{-\frac{\langle f,Cf\rangle}{2}}=\int_{\CS'} e^{i\phi(f)} d\mu_C(\phi),\] where $C$
is the integral kernel of a positive, continuous, nondegenerate Euclidean covariant
bilinear form $C$ on $\CS\times\CS$. We see that the special case
$S\{f\}=e^{-\|f\|^2/2}$ is obtained when $C(s,t)=\delta(t-s)$, called
the ``ultralocal'' covariance. For each $\phi\in\CS'$, the random
variables $\{\phi(f):f\in\CS\}\subset L^2(\CS',{\cal B},\mu)\}$
form a Gaussian family with mean zero and covariances
\begin{equation} E_\mu[\phi(f)\phi(g)]=\int_{\RR_+\times\RR_+}f(s)g(t)C(s,t)ds dt=\langle
f,g\rangle.
\end{equation}

The theory of martingales makes its appearance in white noise
calculus through the concept of Wick ordered random variables. To
introduce them, recall that the Hermite polynomials are defined by
\[H_n(x)=(-1)^ne^{\frac{x^2}{2}}\frac{d^n}{dx^n}(e^{-\frac{x^2}{2}}),\]
from which it is easy to see that they satisfy the recurrence
relation
\begin{eqnarray*}
\left(x-\frac{d}{dx}\right)H_{n-1}(x)&=&H_n(x),\quad
n=1,2,\ldots\\
H_0(x)=1.
\end{eqnarray*}
The first couple of Hermite polynomials are then
\begin{eqnarray*}
H_0(x)&=&1 \\
H_1(x)&=&x \\
H_2(x)&=&x^2-1\\
H_3(x)&=&x^3-3x, \ldots
\end{eqnarray*}

We begin by defining the {\em Wick ordered} exponential for any
$f\in\CS$ to be
\begin{equation} :e^{\phi(f)}: = e^{\phi(f)-\|f\|^2/2}.
\end{equation}
Then we have
\begin{proposition} For any $f\in\CS$,
\begin{equation}
\label{wickexp} :e^{\phi(f)}: =\sum_{n\ge 0}\frac{\|f\|^n}{n!}
H_n\left(\frac{\phi(f)}{\|f\|}\right),
\end{equation} where $H_n$ is the $n$th Hermite polynomial. Moreover, for any $f,g\in
\CS$,
\begin{equation}
\label{Horthog}
E_\mu\left[H_n\left(\frac{\phi(f)}{\|f\|}\right)H_m\left(\frac{\phi(g)}{\|g\|}\right)\right]=\delta_{nm}\
(n!) \left(\frac{\langle f,g\rangle}{\|f\|\  \|g\|}\right)^n
\end{equation}
\end{proposition}

\proof For any $a,b\in \RR$ we have the absolutely convergent
expansions
\begin{eqnarray} e^{a-b^2/2}&=& e^{-(b-a/b)^2/2+a^2/(2b^2)}\nonumber\\ &=& \sum_{n\ge 0}
\frac{b^n}{n!}(-1)^ne^{a^2/(2b^2)}\frac{d^n}{dx^n}\left(e^{-x^2/2}\right)\Big|_{x=a/b}
\nonumber\\ &=&\sum_{n\ge 0}\frac{b^n}{n!}H_n(a/b)
\end{eqnarray} where the last line makes use of the defining property of the Hermite
polynomials. Using this with $a=\phi(f),b=\|f\|$ gives
\eqref{wickexp}. The orthogonality relation (\ref{Horthog})
follows by expanding the identity
\begin{equation}
 E_\mu[:e^{\phi(f)}::e^{\phi(g)}:]=e^{\langle f,g\rangle}
\end{equation} in powers of $f,g$ and comparing to the expansion derived from
(\ref{wickexp}).

\qed

From this we define the {\em Wick ordered} monomials as the random
variables
\begin{equation}
:\phi(f)^n:=\|f\|^n H_n\left(\frac{\phi(f)}{\|f\|}\right),
\label{wickmon}
\end{equation}
so that linearity and convergent
power series imply that
\[:e^{\phi(f)}:=\sum_{n\geq 0}\frac{:\phi(f)^n:}{n!}.\]
Formally, we express the Wick ordered monomials as
\[:\phi(f)^n:=\int_{\RR_+^n}f(s_1)\dots f(s_n)\ :\phi_{s_1}\dots
\phi_{s_n}:\ ds_1\dots ds_n\] and from the orthogonalization
\eqref{Horthog} we can define the {\em Wick products}
\begin{eqnarray}
:\phi(f_1)\dots\phi(f_n): &=&  \int_{\RR_+^n}f_1(s_1)\dots f_n(s_n)\ :\phi_{s_1}\dots
\phi_{s_n}:\ ds_1\dots ds_n, 
\end{eqnarray}

 $\CH_n$ is defined to be the span of
$\{:\phi(f)^n:|f\in\CS\}$ and consists of precisely the random
variables $(n!)^{-1}\int_{\RR_+^n} \tilde f(s_1,\dots,s_n)
:\phi_{s_1}\dots \phi_{s_n}:\ ds_1\dots ds_n$ where $\tilde f$ lies in
$\tilde T_n$ the $L^2$ completion of the space
 of {\em  symmetric} functions  in $\CS^{\otimes n}$. These   subspaces   form an
orthogonal decomposition of
$L^2(\CS',\CB,\mu)$. In quantum field theory this is known as the Fock space
decomposition of the Hilbert space of quantum states into $n$ particle sectors for
$n\ge 0$ .

The relationship between Wiener measure and the white noise measure is to identify
$\phi(f)\doteq \int f(s) \phi_s ds$ with
$W(f)\doteq \int f(s) dW_s$ for all $f\in \CS\subset L^2(\RR_+)$. One is then lead to
the formal relation $\phi_s=dW_s/ds$, that is ``white noise'' is the ``derivative'' of
Brownian motion.  This identification extends to all orders in the chaos expansion via:
\begin{equation}
 (n!)^{-1}\int_{\RR_+^n}\tilde  f(s_1,\dots,s_n) :\phi_{s_1}\dots
\phi_{s_n}:\ ds_1\dots ds_n=\int_{\Delta_n} f(s_1,\dots,s_n)dW_{s_1}\dots
dW_{s_n}
\end{equation} where the bijection $\tilde f\leftrightarrow  f$ between $\tilde T_n$ and
$L^2(\Delta_n)$ is the restriction map and its inverse. Finally, this leads to the
identification
$L^2(\CS',\CB,\mu)\equiv L^2(\Omega,\CF_\infty,P)$.

One useful consequence of Wick products and the chaos expansion is that conditional
expectations can be handled systematically.

\begin{proposition} Let $f\in L^2(\RR_+)$ and $t\in [0,\infty)$. Then
\begin{eqnarray} E[:e^{\phi(f)}:|\CF_t]&=&:e^{\phi(f_{(t)})}:\nonumber\\
E[:(\phi(f))^n:|\CF_t]&=&:(\phi(f_{(t)}))^n:
\end{eqnarray}  where $[f_{(t)}](s)=\One(s\le t)f(s)$. In other words, the processes
$:e^{\phi(f_{(t)})}:,\\
 :(\phi(f_{(t)}))^n:$ are martingales.
\end{proposition}

\section{Generating functional for chaos coefficients} We derive Theorem
\ref{genfnthm} in the notation of  the white noise calculus. To begin, we
recall the definition of analyticity for a function between complex
Banach spaces.

\begin{definition} Let $f:B\to C$ where
$B,C$ are complex Banach spaces. Then $f$ is analytic at a point $x\in B$
if 
\begin{enumerate}
\item for any finite set $\{x_1,\dots,x_n\}\subset B$ of points the   function 
$F:\mathbb{C}^n\to\mathbb{C}$ 
\begin{equation}
F(\zeta_1,\dots,\zeta_n)=f(x+\sum_i\zeta_ix_i)\end{equation}
 is analytic at $0\in\mathbb{C}^n$;
\item $f$ is continuous at $x$.
\end{enumerate}
\end{definition}

\noindent{\bf Proof of Theorem 2.1} One can check directly that the map
$h\mapsto :e^{\phi(h)}:$ is entire analytic between
$L^2_\mathbb{C}(\RR_+)$ and $L^2_\mathbb{C}(\CS',\CB,\mu)$. Similarly, $Z_X(h)$ is
analytic for any
$h\in L^2_\mathbb{C}(\RR_+)$.

Finally, to verify (\ref{Xexp}), it is enough to take the expectation of
the equation multiplied by  $:\phi(g)^n:$ for arbitrary $g\in
L^2(\Delta),n\ge 0$:
\begin{eqnarray} E[X:\phi(g)^n:] &=&\frac{d^n}{d\lambda^n}Z_X(\lambda
g)\Big|_{\lambda=0}\nonumber\\ &=&n!\int_{\Delta_n}
f^{(n)}_X(s_1,\dots,s_n)g(s_1)\dots g(s_n)ds_1\dots ds_n
\end{eqnarray} whereas by (\ref{Horthog}),
\begin{eqnarray} &&E[\sum_{m\ge 0} (m!)^{-1}\int_{\RR_+^m}
f_X^{(m)}(s_1,\dots,s_m):\phi(s_1)\dots\phi(s_m):ds_1\dots
ds_m:\phi(g)^n:]\nonumber\\ &&\hspace{.5in}=n!\int_{\RR_+^m} g(s_1)\dots
g(s_n)f_X^{(n)}(s_1,\dots,s_n)ds_1\dots ds_n
\end{eqnarray}

\qed

\section{The exponential quadratic formula}

Here we prove the formula (\ref{change}) in the context of
one--dimensional white noise calculus. The proof extends easily to the
multidimensional case.
\medskip

\noindent{\bf Proof of Proposition \ref{measchange}:\ } 
\begin{enumerate}
\item (Step 1) First we note that $Y\in \CC^+$ can be approximated in the Hilbert
space
$\CH_{\le 2}$ by random variables of the form
\begin{equation} Y=A + \sum_{i=1}^M[b_i
\phi(g_i)+c_i:\phi(g_i)^2:/2]
\end{equation} with $\{g_i\}$ a finite orthonormal set in $L^2(\Delta)$
and numbers
$c_i>-1,b_i$. Then $X_i=\phi(g_i)$ form a collection of independent
$N(0,1)$ random variables. For $Y$ of this type, there is a factorization
\begin{equation} E[e^{-Y}]=e^{-A}(2\pi)^{-M/2}\prod_i
\int_\RR\exp[-b_i x-(1+c_i)x^2/2+c_i/2)]dx\end{equation} into one
dimensional Gaussian integrals. Each integral gives the factor
\begin{equation} (2\pi)^{1/2}(1+c_i)^{-1/2}
\exp\frac12[c_i+b_i^2(1+c_i)^{-1}]\end{equation} leading to a formula for
$E[e^{-Y}]$ which agrees with (\ref{measchange}) for $Y$ of this form.
\item (Step 2) Since the formula is true for $Y$ in a dense subset of $\CC^+$, it is
now enough to prove that the map $Y\mapsto E[e^{-Y}]$ is continuous in $C$
provided the kernel $C$ satisfies the stated conditions. By  the definition of
the Carleman--Fredholm determinant 
\begin{equation}
\mbox{det}_2(1+C)=\exp[\mbox{Tr}[\log(1+C)-C].
\end{equation}
Since $\log(1+x)-x=\CO (x^2)$ for $x\to 0$    we see that $\det_2(1+C)$
is well-defined and continuous for $C>-1$ and Hilbert--Schmidt. Therefore, the entire right hand side
of (\ref{measchange}) is continuous in $C$ for $Y\in\CC^+$.
\end{enumerate}

\qed

\vspace{0.2in}
\noindent
{\bf Acknowledgements:} We thank L.P. Hughston, D. Brody and L.C.G. Rogers for  discussions and suggestions. We are grateful
for the comments by the participants of the Financial and Insurance Mathematics talks at ETH Zurich, May 2003, the Probability Seminars at University of Cambridge, June 2003, and the Financial Mathematics and Applied Probability seminars at King's College London, July 2003, where
this work was presented.


\begin{thebibliography}{99}

\bibitem{BroHug03}
D.~Brody and L.~P.~Hughston.
\newblock Chaos and coherence: a new framework for interest rate modelling.
\newblock  Submitted to {\em Proc. Roy. Soc.}, 2003.

\bibitem{CoInRo85b}
J.~C. Cox, J.~E. Ingersoll, Jr., and S.~A. Ross.
\newblock A theory of the term structure of interest rates.
\newblock {\em Econometrica}, 53(2):385--407, 1985.

\bibitem{FleHug96}
B.~Flesaker and L.~P.~Hughston.
\newblock Positive interest.
\newblock {\em Risk Magazine}, 9(1):46--49, 1996.

\bibitem{GliJaf81}
J.~Glimm and A.~Jaffe.
\newblock {\em Quantum physics}.
\newblock Springer-Verlag, New York, 1981.
\newblock A functional integral point of view.

\bibitem{HugRaf03}
L.~P.~Hughston and A.~Rafailidis.
\newblock A chaotic approach to interest rate modelling.
\newblock Preprint, 2003.

\bibitem{Jams95}
F.~Jamshidian.
\newblock A simple class of square-root interest rate models.
\newblock 2:61--72, 1995.

\bibitem{JinGla01}
Y.~Jin and P.~Glasserman. 
\newblock Equilibrium positive interest rates: a unified view.
\newblock {\em Rev. Fin. Studies}, 14(1):187--1214, 2001.

\bibitem{Magh96}
Y.~Maghsoodi.
\newblock Solution of the extended {CIR} Term Structure and Bond Option Valuation.
\newblock {\em Math. Finance}, 6:89-109, 1996. 

\bibitem{Nual95}
D.~Nualart.
\newblock {\em The {M}alliavin calculus and related topics}.
\newblock Probability and its Applications. Springer-Verlag, New York, 1995.

\bibitem{Okse96}
B.~{\O}ksendal.
\newblock An introduction to {M}alliavin calculus with applications to
  economics.
\newblock Unpublished lecture notes available at
  \verb+http://www.nhh.no/for/dp/1996/index.htm+, 1996.

\bibitem{Roge97}
L.~C.~G. Rogers.
\newblock The potential approach to the term structure of interest rates and
  foreign exchange rates.
\newblock {\em Math. Finance}, 7(2):157--176, 1997.

\bibitem{Rutk97}
M.~Rutkowski.
\newblock A note on the {F}lesaker--{H}ughston model of the term structure of interest rates.
\newblock {\em App. Math. Finance}, 4:151--163, 1997.

\bibitem{Yor92a}
M. Yor.
\newblock{\em Some aspects of {B}rownian motion}.
\newblock Lectures in Mathematics ETH Zurich. Birkh{\"a}user--Verlag, Basel, 1992.

\end{thebibliography}
\end{document}